# OPTIMAL DOUBLING THRESHOLDS IN BACKGAMMON-LIKE STOCHASTIC GAMES

HAORU JU, DANIEL LEIFER, STEVEN J. MILLER, SOORAJ A. PADMANABHAN, CHENYANG SUN, LUKE TICHI, BENJAMIN TOCHER, AND KILEY WALLACE

*It is a pleasure to dedicate this paper to Backgammon Grandmaster Art Benjamin, whose kindness and insights all authors on this paper have benefited from over the years (directly and indirectly).*

ABSTRACT. We study variants of a stochastic game inspired by backgammon where players may propose to double the stake, with the game state dictated by a one-dimensional random walk. Our variants allow for different numbers of proposals and different multipliers to the stake. We determine the optimal game state for proposing and accepting, giving analytic solutions in many variants. We also introduce a 3-player generalization of the game and prove basic results about its behavior, in addition to providing a simulation.

## CONTENTS

## 1. INTRODUCTION

1.1. **Background on Backgammon.** Backgammon is a popular board game where two players bet on a stake by alternately rolling dice to dictate the movement of counters. Each player is assigned a collection of counters, with the objective being to move them all into a designated region before the opponent accomplishes the same. At the beginning of each turn, a player rolls 2 dice and moves the counters under constraints that depend on the outcomes of the dice. Unlike a deterministic game such as chess, a great part of backgammon strategy involves taking calculated risks in the face of uncertainty.

Aside from the counters and the dice, a doubling cube introduces interesting dynamics to the game and adds a layer of strategy beyond moving the counters. At any point in time, a player may propose to double the current stake of the game. Either player could propose the first double; after that point, a player could only propose a double after the opponent does, i.e., the double proposals must alternate between the players. Once proposed, the opponent has the choice to accept the double–confirming the doubling of the stake and allowing the game to continue–or to decline the double (to "fold"), which forfeits the game and the current stake.

An interesting question naturally arises: when should one propose a double, or accept a double proposed by the opponent? It is clearly disadvantageous to propose a double while losing, so it might make sense to propose double when winning by a great amount. However, if a double is proposed when victory is imminent, the opponent may elect to fold, nullifying the advantage of the proposed double but saving the

---

*Date*: October 16, 2024.

proponent from a small chance of loss. We thus expect the optimal moment to double to be a compromise between the two influences.

1.2. **Stochastic game as a model of Backgammon.**

1.2.1. *The game.* The games we study are based on one introduced by Keeler and Spencer in [KS], which is inspired by backgammon. [1]

In this stochastic game, there are no pieces to move. The *game state*, as a function of game time, is modeled by a Brownian motion $p = p(t)$ with no drift that begins at $0.5$ and stays in $[0, 1]$, where reaching $0$ indicates loss for Player A, and reaching $1$ indicates victory. The initial stake is normalized to be $1$, and the only actions the players may take are to propose or accept/decline a double. To avoid technical issues, we allow doubles to be offered even at the extreme states $0$ and $1$ (and the payoff made afterward), which makes little practical difference from offering at $0 + \epsilon$ and $1 - \epsilon$ for arbitrarily small $\epsilon$. Thus, a double may be proposed at any game state in the closed interval $[0, 1]$.

After a double is proposed, the right to propose the next double goes to the opponent. At this point, the opponent may choose to fold, forfeiting the current stake and the game, or to accept the double and play on with twice the current stake. The game ends when a double is declined, or when the game state reaches $0$ or $1$ and neither player would double.

1.2.2. *Game state and evaluation.* It is convenient to denote the game state $p$ instead by an ordered pair $(s_A, s_B) := (p, 1 - p)$, where $s_A$ is Player A's game state, and $s_B$ Player B's game state, collectively referred to as *subjective game states*. This is to notationally simplify the analysis by eliminating the need to refer to a specific player: the statement "Player A should double at game state $p$, and Player B should double at game state $1 - p$" simplifies to "a player should double at (subjective) game state $p$". From this point onward, all references to game states are subjective unless indicated otherwise.

A useful quantity to study is the *game evaluation* $(p_A, p_B)$, where $p_A$ (resp. $p_B$) is the probability that the upcoming random walk reaches $1$ before $0$ (resp. $0$ before $1$). The evaluation can be thought of as the probability of winning for each player, if the game had no doubling. We may refer to a player's evaluation, taken to mean the corresponding component of the evaluation tuple, i.e., Player B's evaluation being $p_B$.

As it turns out, the game state and evaluation under these rules are equivalent: $p_A = s_A$, and $p_B = s_B$. This is an immediate consequence of the consistency property of one-dimensional Brownian motion: within any interval, the probability of reaching the right endpoint before the left equals the proportion of the interval to the left of the current location. This consistency property is proved in [Don], though a discrete version is stated and proved in Section 4, as Lemma 4.2.

We seek to determine game states–equivalently, evaluations–where it is most advantageous to propose or accept/decline a double, under different, more general game rules.

1.2.3. *Maximization of expected payoff and a principle of indifference.* We assume that all players seek to maximize their expected payoff in the face of uncertainty. This assumption on player behavior has important implications for doubling. Since the expected payoff for a player is an increasing function of their game evaluation, there exists an evaluation threshold such that a player should offer a double exactly when that threshold has been surpassed. Doubling before that threshold offers the opponent too much chance for a comeback, and waiting to double far past that threshold wastes an opportunity to end the game immediately and eliminate the risk of loss, since a rational opponent would prefer, under a severe disadvantage, to forfeit the current stake by folding than to accept the double and risk losing twice the stake. In [KS] the following is proved; a proof is also included in Appendix A.

---

[1] The actual game of backgammon is more complex than this game, but it serves as a useful model. While the strategy of moving counters and randomness of dice rolling have been simplified greatly and replaced with Brownian motion of the game state, the decision of when to double and to accept/fold remains.

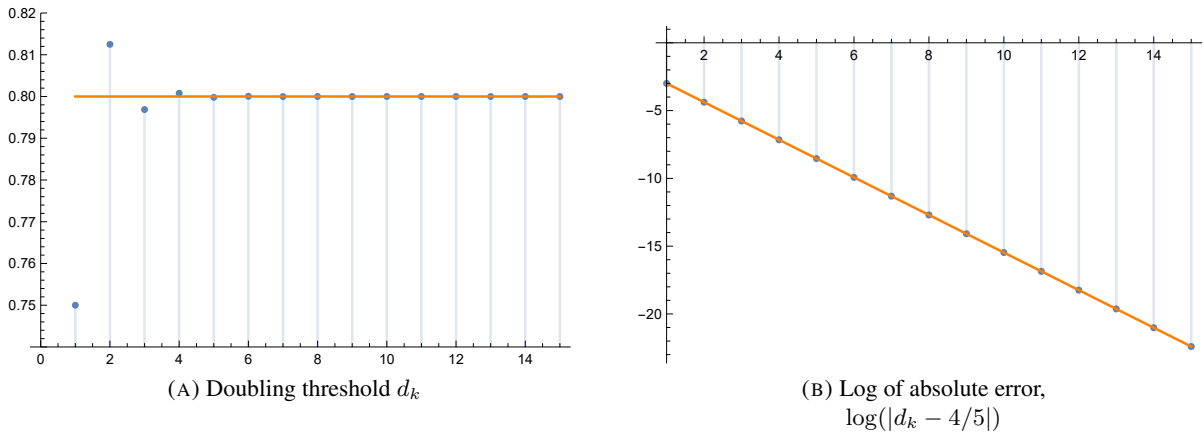

(A) Doubling threshold $d_k$

(B) Log of absolute error, $\log(|d_k - 4/5|)$

FIGURE 1. The first plot is of double-proposal threshold $d_k$ with $k$ doubles remaining. Also plotted is the limit (orange), which equals the threshold for indefinitely many doubles remaining, studied in [KS]. The second plot is the log of the absolute error from the limit $4/5$, i.e., $\log(|d_k - 4/5|)$, with the orange line being $\log(1/5) + x\log(1/4)$.

**Theorem 1.1.** *Let $\alpha$ be Player A's optimal game state for proposing a double, and $\beta$ be Player A's game state where Player B is indifferent between accepting and declining a double. Then $\alpha = \beta$.*

This result can be used for finding the optimal state for doubling by instead computing the point at which the opponent is indifferent between accepting and folding; this indifference can be used to set up an equality from which one can solve for the optimal state.

The analysis in [KS] adopts the rule that either player could propose indefinitely many doubles, and concludes that $0.8$ is the best evaluation to propose a double. We extend this analysis to variants with different rules for modifying the stake.

1.3. **Main results.** In Section 2, we study variants of the stochastic game where the number of doubles remaining is changed to some arbitrary positive integer. The analysis extends that of [KS], and culminates in a closed-form expression for the optimal doubling point given an arbitrary number of doubles remaining.

**Theorem 1.2.** *Suppose there is a total of $k \geq 1$ potential doubles remaining in the game. The optimal game state for doubling is given by*

$$d_k = \frac{4}{5} + \frac{1}{5}\left(-\frac{1}{4}\right)^k. \quad (1.1)$$

Figure 1 plots the optimal game state, as a function of the number of doubles remaining, against its limiting value of $0.8$. Note that the thresholds for an odd number of doubles remaining stay below the limiting threshold, while the even thresholds stay above. This indicates that having the final double is an advantage that grows as fewer doubles remain.

In Section 3, we study variants where the stake changes by a multiplier other than 2 when the "doubling cube" is played. We introduce the $(x, y)$-cube parameterized by real numbers $x, y$: when the cube is played, upon acceptance the current stake is multiplied by $y$, and upon folding the opponent forfeits $x$ times the current stake. The default doubling cube is the special case $(x, y) = (1, 2)$.

We determine the optimal doubling strategy for all $(x, y)$ with $x \in \mathbb{R}$, $y > 0$. We first treat some regions of $(x, y)$-space that yield trivial decisions; for the remaining regions, we prove the following principle of indifference, which generalizes Theorem 1.1.

**Theorem 1.3.** *Let $(x, y)$ be an arbitrary pair of real numbers with $y > 0$. Let $\alpha$ be Player A's optimal evaluation for proposing a double, and let $\beta$ be Player A's evaluation threshold where Player B is indifferent between accepting and folding. If $(x, y)$ falls in one of the regions*

- $R_+ := \{(x, y) \in \mathbb{R}^2 \mid 0 < x < y,\ y > 1\}$,
- $R_- := \{(x, y) \in \mathbb{R}^2 \mid 0 < -x < y,\ 0 < y < 1\}$,

*then $\alpha = \beta$.*

Both regions are illustrated in Figure 4. If $(x, y) \in R_+$, we call the doubling cube a *raising cube*, and in the case of $R_-$, a *reducing cube*, while the term *doubling cube* serves as an umbrella term for both. An offer of the doubling cube, a *double*, specifies to a *raise* and a *reduction*. A raising cube can be used by the winning player to enlarge the payoff, while a reducing cube can be used by the losing player to cut the loss. Note that the raising region $R_+$ contains the pair $(1, 2)$ corresponding to the standard doubling cube.

In the spirit of Theorem 1.2, we have a closed-form expression for the optimal state to double, whether the cube is a raising cube or a reducing cube.

**Theorem 1.4.** *Suppose $(x, y) \in R_+$. The optimal game state $d_k$ to raise, given there are $k \geq 1$ proposals left to play, is*

$$d_k = \frac{y(x+1)}{2y + xy - x} + \frac{x(1+x)(y-1)}{2(2y + xy - x)} \left(-\frac{y-x}{y(1+x)}\right)^k. \tag{1.2}$$

*Furthermore, raising at $d_k$ is preferable to not raising at all.*

Note that the second part of the theorem is not trivially true: if $x < 1$, it might be preferable to not raise when victory is imminent, in order to collect the full stake rather than the reduced stake $x$ if the opponent folds.

In analogy, the following results holds for the reducing cube.

**Theorem 1.5.** *Suppose $(x, y) \in R_-$. The optimal game state $d_k$ to reduce, given there are $k \geq 1$ reductions left to play, is*

$$d_k = \frac{y + x}{2y - xy + x} + \frac{(-x)(1-x)(1-y)}{2(2y - xy + x)} \left(-\frac{y+x}{y(1-x)}\right)^k. \tag{1.3}$$

*Furthermore, reducing at $d_k$ is preferable to not reducing at all.*

Figure 2 gives a contour plot of optimal raising points $d_k$ in (1.2) over $R_+$, for $k = 1,\ 2,\ \infty$. The plots approach $1/2$ along the boundary $x = 0$, and $1$ along the boundary $y = x$. In parallel, Figure 3 plots the optimal reducing points $d_k$ in (1.3). The plots approach $0$ along the boundary $y = -x$, and $1/2$ along the boundary $x = 0$, with a discontinuity at $(0, 0)$.

In Section 4 we introduce a three-player variant of the game based on a discrete random walk in a equilateral triangle $\triangle ABC$, where a player wins by reaching the corresponding vertex. The game state is a point in the triangle, expressible in barycentric coordinates with the vertices being the basis. We show that at any point $P$, the game evaluation for an arbitrary Player A, i.e., the probability that the random walk from $P$ first reaches vertex $A$, is given by the barycentric $A$-component of $P$. This also serves as a two-dimensional discrete analogue to the consistency property of 1-dimensional Brownian motion.

We then extend 2-player doubling rules, and prove a principle of indifference (Theorem 4.3) analogous to Theorem 1.1. We use it to compute the optimal doubling points and acceptance thresholds for a game with one double allowed, observing that the thresholds may depend on the evaluation of all players, not just one's own. We finish by proposing alternate models for a 3-player game and presenting a computer simulation of the triangle random walk.

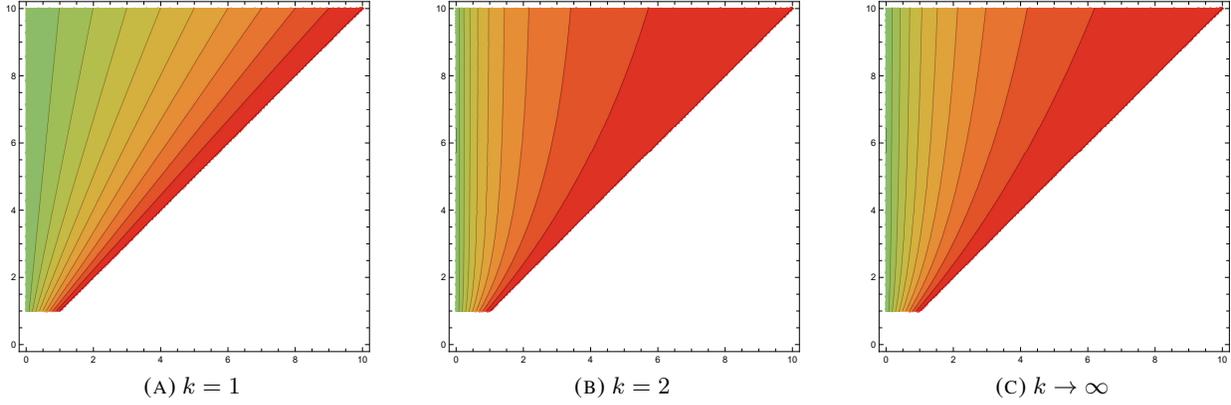

FIGURE 2. Contour plots over a portion of the raising region $(x, y) \in R_+ \cap [0, 10]^2$ of $d_k$ given in (1.2) for $k = 1$, $k = 2$, and the limiting value as $k \to \infty$. The horizontal axis is $x$, and the contours are spaced by $1/20$, from $1/2$ (green) to $1$ (red).

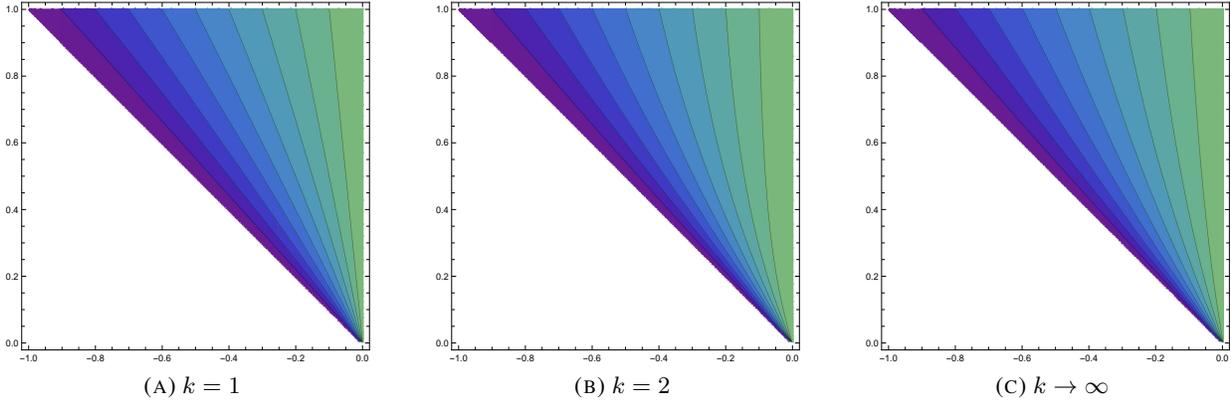

FIGURE 3. Contour plots over the reducing region $(x, y) \in R_-$ of $d_k$ given in (1.3) for $k = 1$, $k = 2$, and the limiting value as $k \to \infty$. The horizontal axis is $x$, and the contours are spaced by $1/20$, from $0$ (purple) to $1/2$ (green). Note the discontinuity at $(0, 0)$.

## 2. TWO-PLAYER GAMES WITH VARYING NUMBER OF DOUBLES

**2.1. One double allowed in total.** When only one double is allowed, let $p$ be a player's optimal game state for doubling. By Theorem 1.1, this occurs when the expected payoff is the same whether the double is accepted or not.

Normalizing the initial stake to be 1, the expected payoff in the case of an accepted double is $2p - 2(1-p)$, while the expected payoff in the case of a declined double is 1. Setting $2p - 2(1 - p) = 1$ yields $p = 3/4$.

**2.2. Up to $n$ alternating doubles remaining.** Here we treat a more general set of rules, and suppose that there are a total of $n$ doubles allowed in a game. We allow the game to possibly end before all $n$ doubles are played.

Let $a_k$ denote the evaluation threshold for a player to accept/reject the $k^{\text{th}}$-to-last double, and $d_k$ the evaluation threshold for a player to propose the $k^{\text{th}}$-to-last double. We wish to determine $a_k, d_k$.

*Proof of Theorem 1.2.* At evaluation $a_k$, a player is indifferent between accepting and declining a potential $k^{\text{th}}$ double proposed by the opponent; either decision results in the same expected payoff. Declining results in a payoff of $-1$, i.e., losing the initial stake. Upon accepting the double, one of two items may happen: losing (with payoff $-2$) without being able to double again, or recovering to some critical level $d_{k-1}$ where the player can propose the $(k-1)^{\text{st}}$-to-last double (as long as $k \geq 2$, and the opponent is indifferent between accepting and folding. By indifference, the payoff is 2 in this case.

Given the current evaluation is $a_k$, the probability of reaching $d_{k-1}$ before 0 is $(a_k - 0)/(d_{k-1} - 0)$ by the consistency property of Brownian motion, i.e., the probability of reaching the right endpoint of an interval before the left endpoint is the proportion of the interval traveled from the left endpoint. Correspondingly, the probability of reaching 0 first is $(d_{k-1} - a_k)/d_{k-1}$. The expected payoff for accepting the double is thus $(-2)((d_{k-1} - a_k)/d_{k-1}) + 2(a_k/d_{k-1})$. Setting this equal to the expected payoff for folding, we have

$$-1 \;=\; (-2)\frac{d_{k-1} - a_k}{d_{k-1}} + 2\frac{a_k}{d_{k-1}}, \tag{2.1}$$

In the case $k = 1$, a player could not later double back after recovery; the evaluation must go to 1 in order to win, since they could not double at a lower evaluation to force the opponent to fold. We may encode this necessity by setting $d_0 = 1$.

Now we consider the threshold for proposing a double. By Theorem 1.1, the threshold for proposing a double corresponds to the opponent's threshold for accepting, i.e., $d_k = 1 - a_k$. Substituting this into (2.1) gives

$$-1 \;=\; -2\left(\frac{d_{k-1} - (1 - d_k)}{d_{k-1}}\right) + 2\left(\frac{1 - d_k}{d_{k-1}}\right), \tag{2.2}$$

which simplifies to

$$d_k \;=\; 1 - \frac{d_{k-1}}{4}. \tag{2.3}$$

The initial condition $d_0 = 1$ thus allows the explicit determination of $d_k$ and $a_k$ for all $k$.

Observe that this recurrence stabilizes if one were to take $4/5$ instead as the initial condition, which makes one suspect that $4/5$ might be the limiting value. Notice that $4/5$ is the established threshold for doubling if the number of doubles are unlimited. We could try to keep track of the error term $\epsilon_k := d_k - 4/5$, and substitute the recurrence (2.3) to obtain

$$\begin{aligned}
\epsilon_k &= 1 - \frac{d_{k-1}}{4} - \frac{4}{5} \\
&= \frac{1}{5} - \frac{d_{k-1} - 4/5 + 4/5}{4} \\
&= \frac{1}{5} - \frac{\epsilon_{k-1}}{4} - \frac{1}{5} \\
&= -\frac{\epsilon_{k-1}}{4},
\end{aligned} \tag{2.4}$$

and we see that the error term satisfies $\epsilon_k = (-1/4)^k \epsilon_0$, undergoing a simple exponential decay. Setting $\epsilon_0 = 1 - 4/5 = 1/5$ thus gives

$$\begin{aligned}
d_k &= \frac{4}{5} + \epsilon_k \\
&= \frac{4}{5} + \frac{1}{5}\left(-\frac{1}{4}\right)^k.
\end{aligned} \tag{2.5}$$

$\square$

## 3. Two-player games with varying stakes

We now explore more general game rules that allow for the stakes to vary. There are two parameters at hand: the multiplier $y$ to the stake upon accepting a double, and the multiplier $x$ to the forfeited stake upon folding. The results of this section mirror those of Section 2.

### 3.1. Trivial regions.
We first consider values of $(x, y)$ that result in trivial decisions.

**Theorem 3.1.** *Let $x \in \mathbb{R}$ be the multiplier to the forfeited stake upon declining a double, and $y \in \mathbb{R}^+$ be the multiplier to the stake upon accepting a double. For certain regions in the $(x, y)$-space, the following optimal doubling strategy holds.*

*(1) If $x \geq y \geq 0$ and $y > 1$, one should double at evaluation $p = 1$.*
*(2) If $x \geq 0$ and $y \leq 1$, one should double at evaluation $p = 0$.*
*(3) If $x < 0$, and $-x > 1$ or $y > 1$, then one should never double.*
*(4) If $-1 \leq x \leq 0$ and $y \leq -x$, one should double at evaluation $p = 0$.*

Figure 4 illustrates the corresponding regions and summarizes the results for the optimal time to double.

*Proof of Theorem 3.1.* Suppose that $x \geq y \geq 0$. If offered a double, a player always fares better to accept scaling the stake by $y$ than to fold, which forfeits at least as much with certainty. If $y > 1$, then the double should be offered at 1. At $p = 1 - \epsilon$, the evaluation could eventually go to 1 or 0: in the former case one could have waited till $p = 1$ to double; in the latter case, one has introduced unnecessary risk, since the opponent would prefer to accept and stay in the game. Thus, waiting till $p = 1$ is strictly better. If $y \leq 1$, a similar argument shows that the double should be offered at 0, in order to not forfeit a chance of coming back. It is easy to check that in both cases, offering a double at an extreme point is better than not offering, which gives up an opportunity to scale the payoff favorably.

Suppose that $x \geq 0$ and $y \leq 1$. If offered a double, accepting results in reducing the stake, and folding forfeits some nonnegative payoff. At game state $p > 1/2$, a player should not double since the opponent may choose to accept and reduce the stake. If a player doubles at game state $p \leq 1/2$, the opponent would always accept to avoid forfeiture. If $p > 0$, there is a nonzero probability that the doubler might make a comeback and regret doubling and reducing the stake; in the case of not making a comeback, the player could double at $p = 0$ to reduce the stake. Thus, $p = 0$ is the optimal time to double.

Suppose that $-x \geq y \geq 0$. If a double is offered, folding immediately pays more than the post-acceptance stake, so it is always advantageous to fold. If $x \leq -1$, then folding pays at least as much as the current stake, so a double should never be offered. If $x > -1$, then a double should be proposed only if $p = 0$, since doubling at any other time forfeits a nonzero probability of making a comeback. Conversely, at $p = 0$ it is advantageous to double and reduce the amount forfeited when the opponent folds.

Suppose that $-x \geq 0$ and $y \geq 1$. If offered a double, accepting it enlarges the current stake, and folding wins some nonnegative payoff. At game state $p \geq 1/2$, a player should not offer a double, since the opponent may simply fold and receive a nonnegative payoff, rather than the negative expected payoff pre-doubling. At evaluation $p < 1/2$, one should not double either since the opponent may choose to accept and enlarge the stake, which is worse than not doubling since the doubler has a disadvantageous game state. Thus, a player should never double in this case regardless of game state. □

### 3.2. Nontrivial regions.
For the remaining nontrivial values of $(x, y)$–that is, those in one of the regions

- $R_+ := \{(x, y) \in \mathbb{R}^2 \mid 0 < x < y, \ y > 1\}$,
- $R_- := \{(x, y) \in \mathbb{R}^2 \mid 0 < -x < y, \ 0 < y < 1\}$,

Theorem 1.3, a principle of indifference generalizing Theorem 1.1, is proven below.

*Proof of Theorem 1.3.* Let $p$ be Player A's game state, and suppose that Player A proposes a double. Let $E_B(p)$ be the expected payoff for Player B upon accepting the double. By our game rules $E_B(p)$ is a decreasing function of $p$ with $E_B(0) = y$ and $E_B(1) = -y$, which corresponds to winning/losing the updated stake. To decide whether to accept the double, Player B checks $E_B(p)$ against the rejection payoff $-x$; since $|x| < y$ in both regions $R_+$, $R_-$, there exists a threshold $\beta$ where $E_B(\beta) = -x$, and Player B is indifferent between acceptance and rejection.

Case 1: $(x, y) \in R_+$. At any point $p$, if Player B prefers folding to accepting (i.e., $p > \beta$), then $p$ is not the optimal time to double. If $x \geq 1$, then Player A could have doubled at an earlier time when Player B would fold, which is the best outcome for Player A since Player B can guarantee losing $x$ or less. If $x < 1$, then Player A should wait, since there is a nonzero probability that $p$ reaches 1 before $\beta$, which results in the full payoff of 1 instead of $x$.

Neither is any $p < \beta$ the optimal time to double, since waiting till $p \geq \beta$ eliminates Player B's chance of coming back after accepting. Thus $\alpha = \beta$.

Case 2: $(x, y) \in R_-$. At any point $p$, if Player B prefers to fold and receive $-x > 0$, i.e., $E_B(p) < -x$, then Player $A$ would have fared better to accept and keep playing for a chance of comeback. If Player B prefers to accept and continue the game, i.e., $E_B(p) > -x$, then it would have been better for Player A to double earlier and cut the losses to $-x$. Thus the optimal point of doubling coincides with the opponent's point of indifference between accepting and folding, and $\alpha = \beta$. □

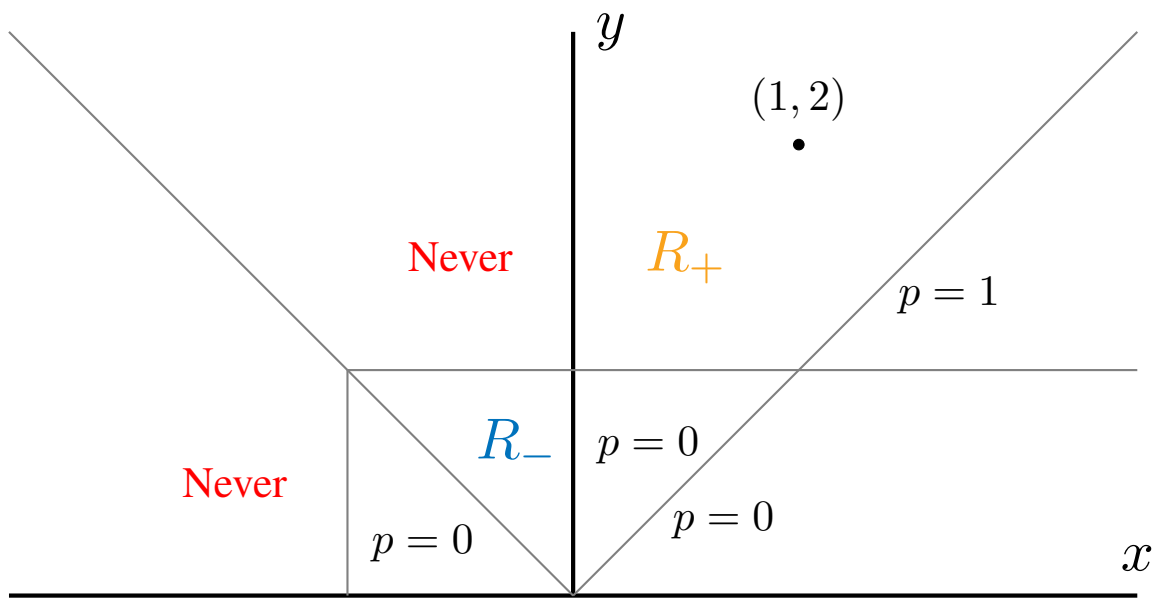

FIGURE 4. Regions in $(x, y)$-space labeled with the optimal time to double in the trivial cases. The thick lines are $x, y$-axes; the diagonal lines are $y = \pm x$, and the thin horizontal line is $y = 1$. The nontrivial regions $R_+$, $R_-$ are as in Theorem 1.3. The point $(1, 2)$ corresponds to the regular doubling cube.

3.2.1. *Raising cube.* We compute the optimal doubling points for $(x, y) \in R_+$.

For a base case, we determine the optimal doubling point for only one double, which includes both raising and reducing.

**Lemma 3.2.** *Suppose $(x, y) \in R_+ \cup R_-$. Suppose that if one doubles at game state $p$, the opponent is indifferent between accepting and folding. Then*

$$p = \frac{1}{2} + \frac{x}{2y}. \tag{3.1}$$

*Proof.* Let $p$ be the the point of indifference as described. In the case of acceptance, the doubler's expected payoff, assuming unit initial stake, is $yp - y(1-p) = y(2p-1)$. Setting $y(2p-1) = x$, the payoff in the case of folding, yields

$$p = \frac{1}{2} + \frac{x}{2y} \tag{3.2}$$

as claimed. □

Note that substituting $y = 2$ and $x = 1$ gives $p = 3/4$, which matches the result for one single offer of the doubling cube.

With the $k = 1$ case prepared, we can now prove Theorem 1.4.

*Proof of Theorem 1.4.* The analysis will be similar to, but slightly more nuanced than that of Theorem 1.2. The complication comes from the possibility that for some values of $(x, y)$ and $k$ doubles remaining, it could be better to never double than to double at the optimal point. We will rule out this possibility at the end of the proof.

For simplicity, we first determine the doubling points for a pair of opponents who do not perform the comparison against never doubling, and always double at the optimal point.

**Lemma 3.3.** *Consider a $(x, y)$-cube game with $(x, y) \in R_+$ with $k$ doubles remaining. Suppose that players always double at evaluation $d_k$ where the opponent is indifferent between accepting and folding. Then*

$$d_k = \frac{y(x+1)}{2y + xy - x} + \frac{x(x+1)(y-1)}{2(2y + xy - x)} \left(-\frac{y-x}{y(x+1)}\right)^k. \tag{3.3}$$

.

The proof, being a more tedious version of the proof for Theorem 1.2, is left to Appendix A.

Now we determine when it is better to double at the optimal point than to never double. At state $p$, not doubling results in an expected payoff of $2p - 1$; in order for doubling to be preferable, this must be no greater than the payoff $x$ of doubling at the optimal point. Solving $2p - 1 \leq x$ yields $p \leq (1+x)/2$; whenever $d_k$ is less than $(1+x)/2$, doubling is preferable to doubling. We need to check that $d_k < (1+x)/2$ holds for all values of $(x, y) \in R_+$ and $k \geq 1$.

In (3.3), observe that the constant term and the coefficient of the exponential term are both positive, since $n > 1$. Furthermore, the exponential term has magnitude less than 1 and alternates in sign, which implies that $d_k$, as a function of $k$, is maximized when $k = 2$. It thus suffices to show that $d_2 < (1 + x)/2$. Fortunately, it holds that

$$\frac{1+x}{2} - d_2$$
$$= \frac{1+x}{2} - \frac{y(x+1)}{2y + x(y-1)} - \frac{x(x+1)(y-1)}{2(2y + x(y-1))} \left(-\frac{y-x}{y(x+1)}\right)^2$$
$$= \frac{x^2(y^2 - 1)}{2y^2(x+1)}, \tag{3.4}$$

which is positive since $y > 1$, completing the proof. □

3.2.2. *Reducing cube.* We compute the optimal reducing strategy for $(x, y) \in R_-$. With the $k = 1$ case prepared, we can now prove Theorem 1.5.

*Proof of Theorem 1.5.* As in the proof of Theorem 1.4, we compute the optimal reducing states assuming that each player would reduce whenever the opponent is indifferent, and then show that doing so is always preferable to not reducing.

**Lemma 3.4.** *Consider a $(x, y)$-cube game with $(x, y) \in R_-$ with $k$ reductions remaining. Suppose that players always reduce at state $d_k$ where the opponent is indifferent between accepting and folding. Then*

$$d_k = \frac{y + x}{2y - yx + x} + \frac{(-x)(1 - x)(1 - y)}{2(2y - yx + x)} \left( -\frac{y + x}{y(1 - x)} \right)^k. \tag{3.5}$$

The proof is similar to the proof of Lemma 3.3 and left to Appendix A.

Now we show that it is better to reduce at the optimal point than to never reduce. At state $p$, not doubling results in an expected payoff of $2p - 1$; in order for doubling to be preferable, this must be no greater than the payoff $x$ of doubling at the optimal point. Solving $2p - 1 \leq x$ yields $p \leq (1 + x)/2$; whenever $d_k$ is less than $(1 + x)/2$, doubling is preferable to doubling. We need to check that $d_k < (1 + x)/2$ holds for all values of $(x, y) \in R_-$ and $k \geq 1$.

In (3.5), observe that the coefficient of the exponential term is positive, since each factor in the numerator is positive, and the denominator agrees in sign with $2y - yx + x = y(1-x) + (y+x)$, where both summands are positive since $y > |x|$. Furthermore, the exponential term has magnitude less than 1 and alternates in sign, which implies that $d_k$, as a function of $k$, is maximized when $k = 2$. It thus suffices to show that $d_2 < (1 + x)/2$. Fortunately, it holds that

$$\begin{aligned}
&\frac{1 + x}{2} - d_2 \\
&= \frac{1 + x}{2} - \frac{y + x}{2y - yx + x} - \frac{(-x)(1 - x)(1 - y)}{2(2y - yx + x)} \left( -\frac{y + x}{y(1 - x)} \right)^2 \\
&= \frac{x^2(1 - y^2)}{2y^2(1 - x)},
\end{aligned} \tag{3.6}$$

which is positive since $y < 1$ and $x < 0$, completing the proof. $\square$

## 4. Three-player variants

It is natural to generalize the existing game to a variant involving 3 players. In the 2-player game, we used a line segment $[0, 1]$ to model the space of game evaluations; in a 3-player setup, it would be natural to choose a two-dimensional region. There are many choices one can make about the boundary shape of the region and the number of directions a single turn can traverse. In this section we mainly explore a random walk in a triangle that leads to neat properties, give results on doubling thresholds, and present a computer simulation, along with a sketch of other random walk models.

**4.1. Random walk in a Triangle.** One natural choice is to use an equilateral triangle and its interior to represent the space of game states, where the vertices $A, B, C$ correspond to the three players. Each point in the triangle can be written in barycentric coordinates as $xA + yB + zC$, or $(x, y, z)$ in short, where $x, y, z \in [0, 1]$, and $x + y + z = 1$. Each coordinate has a geometric interpretation: if the triangle $ABC$ has area 1, each coordinate of $P = (x, y, z)$ is the area of the triangle subtended by the point $P$, and the vertices associated with the remaining coordinates; see Figure 5.

Once a vertex is reached, the corresponding player is declared victorious; this corresponds to a coordinate where one entry equals 1, and the other coordinates are 0, which would be a win in the 2-player game. Additionally, it is also natural to stipulate additionally that once an edge is reached, the player corresponding to the non-incident vertex is ejected from the game, and the game devolves into a two-player game between the incident vertices. This corresponds to a player's coordinate reaching 0, which would be a loss in the 2-player game.

For simplicity we study a discretized game, with the triangle's interior filled in by an equilateral triangle tiling. The game state is modeled by a vertex of the tiling. A point in the interior of the triangle is allowed to move to one of the six adjacent vertices per turn, with equal probability; a point on the boundary is allowed to move only in 2 directions along the boundary, with equal probability, corresponding to the ejection of the remaining player.

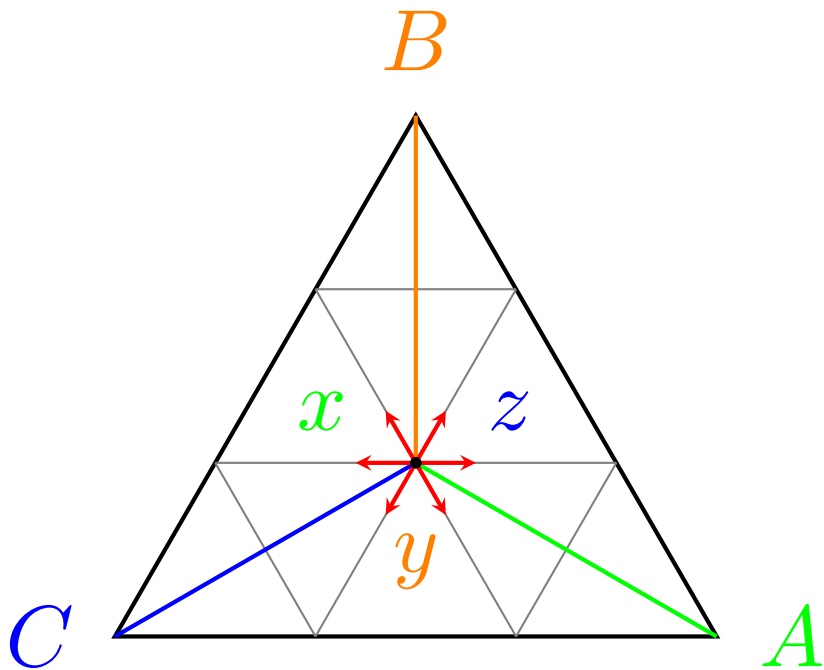

FIGURE 5. A unit-area equilateral triangle with an underlying tiling whose vertices correspond to discretized game evaluations. The coordinates $(x, y, z)$ of an arbitrary point $xA + yB + zC$ correspond to the areas shown. The six arrows correspond to the 6 directions of movement available to an interior point.

Under these game rules, the game evaluation at any point $P$ can be written as $(P_A, P_B, P_C)$, where $P_A$ denotes the probability of the point $P$ reaching vertex $A$ before any other vertex. Clearly $P_A + P_B + P_C = 1$, which reminds one of a barycentric coordinate. We now state and prove a remarkable expression for the game evaluation.

**Theorem 4.1.** *Given a game state $P = xA+yB+zC$, the game evaluation at $P$ is given by $(P_A, P_B, P_C) = (x, y, z)$.*

Before proving the theorem, we first state and prove a discrete analogue to the consistency property of Brownian motion, which is useful toward treating the edges of the triangle.

**Lemma 4.2.** *Let $n$ be a positive integer, and $0 \leq k \leq n$ be an integer. Consider the random walk constructed by adding $+1$ or $-1$ to $k$ independently and with equal probability at each step. Then, the probability that $n$ is reached before $0$ is $k/n$.*

*Proof.* Let $p_k$ ($0 \leq k \leq n$) denote the probability that starting at $k$, $n$ will be reached before $0$. We immediately have $p_0 = 0$ and $p_n = 1$.

For $1 \leq k \leq n-1$, there is an equal chance that $k$ becomes $k-1$ or $k+1$ after the next addition, so $p_k = (p_{k-1} + p_{k+1})/2$. Rearranging gives

$$\begin{aligned} 2p_k &= p_{k-1} + p_{k+1} \\ p_k - p_{k-1} &= p_{k+1} - p_k, \end{aligned} \tag{4.1}$$

which implies the $p_k$'s form an arithmetic sequence due to the common difference. Knowing $p_0 = 0$ and $p_n = 1$ suffices to conclude that $p_k = k/n$. □

Noting that $k/n$ is the barycentric $n$-coordinate of the point $k$ in the basis $\{0, n\}$, we may now prove Theorem 4.1.

*Proof of Theorem 4.1.* We proceed to verify the claim along the boundary of the triangle and then in the interior.

Consider first the vertices $A, B, C$, which correspond to the barycentric coordinates $(1, 0, 0)$, $(0, 1, 0)$, $(0, 0, 1)$, which trivially gives the evaluations of the players since they each represent a game termination state.

Now consider a point $P = (x, y, z)$ in the triangle that is not one of the vertices; we split into two cases depending on the location of $P$. If $P$ is on a boundary, then by the game rules $P$ undergoes a symmetric one-dimensional random walk along an edge, say $\overline{AB}$, i.e., $z = 0$. By Lemma 4.2, the probability $P_A$ of reaching $A$ before $B$ is the barycentric coordinate $x$. Similarly, $P_B = y$, $P_C = z$ by symmetry.

If $P$ is in the interior, then its evaluation $(P_A, P_B, P_C)$ is the average of the evaluations of the 6 incident vertices in the triangle tiling. Observe that this recurrence relation among the triangle vertices is also satisfied by the barycentric coordinates $(x, y, z)$. Since the equality $(P_A, P_B, P_C) = (x, y, z)$ holds on the boundary of the triangle, the recurrence forces the equality to hold everywhere. □

4.2. **The Stake and Doubling Dice.** The introduction of an additional player poses the question of the most natural extensions to the 2-player game rules.

In the 2-player game, a player stands to win or lose by the same amount, which is termed the stake. In the 3-player game, we stipulate that the payment to the winner is split evenly among the remaining opponents. Here we define the stake to be the payment made by each opponent; thus, the winner receives twice the stake.

In the three-player game, a doubling offer, if accepted, doubles the stake of every player. A sufficient condition for the double to be accepted is up to choice; the game can have the double effective only if both opponents accept, or have one acceptance be sufficient and effective. One might argue that it is natural to have one acceptance suffice, since otherwise a losing player may exert undue influence on another player by refusing a double and immediately ending the game.

After the doubling cube is played, the right to the cube rests with the remaining two players. It is not required for both opponents to double before a player could double again.

We prove a principle of indifference for the 3-player game analogous to Theorem 1.1.

**Theorem 4.3.** *A player's optimal time to double occurs when both opponents are first willing to fold.*

*Proof.* Consider the perspective of Player A. Suppose that both opponents are willing to fold if Player A doubles. Since Players B, C are rational, they would not each surrender more than their current stake; thus Player A achieves the optimal outcome (collecting twice the current stake) by doubling. Thus the optimal time for Player A to double occurs when both opponents are first willing to fold, or before; it is not advantageous to wait past this point.

For the other direction, suppose that Player B is willing to accept a double. If Player A doubles at this point, they would experience one of the two possible outcomes: reaching a point when Players B, C would both fold, or losing the game. Player A would have fared better to wait until both opponents would fold to propose the double, which eliminates the risk of losing the game. Thus the optimal time to double agrees with the time when both opponents are first willing to fold. □

4.3. **Decision thresholds with one double.** We determine the decision thresholds for a game that allows for only one double to be played.

**Theorem 4.4.** *In the three-player game with one double allowed, a player should fold if their evaluation is less than $1/6$, and a player should double if both opponents have an evaluation no more than $1/6$.*

*Proof.* Without loss of generality we consider the perspective of Player B, who is indifferent between accepting and folding at evaluation $(x, y, z)$. By folding, if the folding is effective, Player B's payoff is $-1$.

If Player B accepts the double, they may lose via the evaluation reaching $0$, which occurs with probability $1 - y$ with payoff $-2$, or win. Since there are no doubles left, the only way to win is to have the evaluation reach $1$, which happens with probability $y$ and payoff $4$. Thus the overall expected payoff is $4y - 2(1 - y)$. Setting this equal to $-1$ yields $y = 1/6$, which is Player B's accept/reject threshold. The claimed doubling threshold follows from Theorem 4.3. □

**Remark 4.5.** *The game with only one double, as it turns out, is exceptionally easy to study. When more doubles are available, Player B from the proof, after accepting a double, may win without having the evaluation reach 1–it suffices to reach the threshold for proposing a double. As we see in Theorem 4.4, the threshold might not be expressible as Player B's evaluation alone, instead depending on the evaluations of both opponents, forming a curve inside the evaluation triangle. This introduces the implicit problem of determining, for a particle undergoing a random walk, the probability of reaching an arbitrary region before reaching an edge of the triangle.*

Figure 6 illustrates the decision regions obtained in Theorem 4.4.

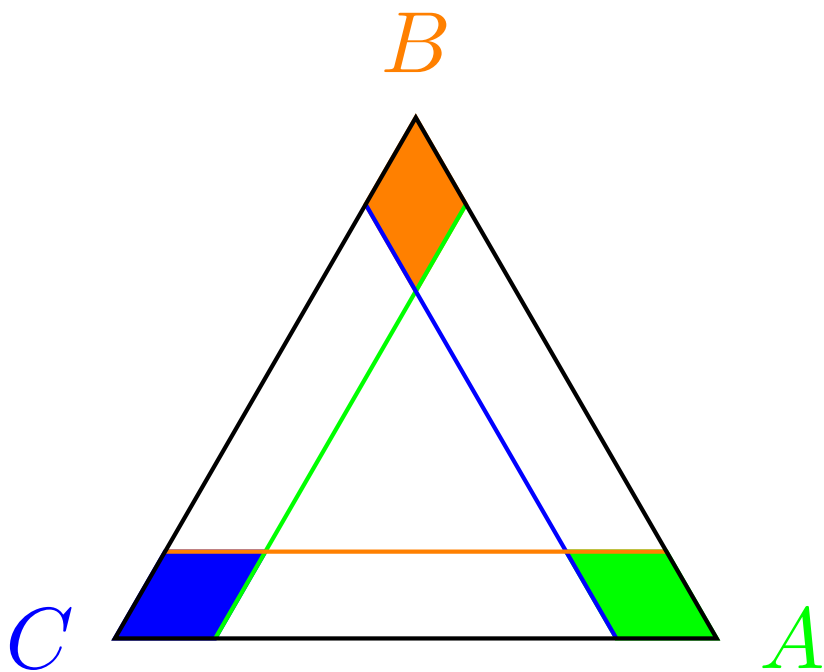

FIGURE 6. The decision regions if one double is allowed, per Theorem 4.4. A player corresponding to a colored vertex should propose a double at the boundary of the filled region, and be willing to fold at the faraway line of the same color.

In Appendix B, we give a computer simulation of the random walk described in this section. We created several functions which allowed us to (i) simulate a round of the 3-player game, (ii) simulate several rounds to find the probability of each player winning from a certain point, (iii) create a function to test different bounds for doubling and explore different winners and winning stakes, and (iv) compare the probability of each player winning at a certain point compared to the area of their triangle.

4.4. **Alternative 3-player models.** Aside from the random walk in a triangle, other models of a 3-player game may be of interest.

Another random-walk model can be constructed using the regular hexagon as the space of game states, with three non-adjacent edges each representing a zone of victory for the three players.

The random walk starts from the geometric center of the hexagon and allows moves in three directions with probability $1/3$ during each discrete time step. These three directions are orthogonally toward the victory zone of each player, see Figure 7. When reaching a non-winning edge, the game designer may choose to eject the antipodal player, or to assign a probability of $1/2$ to each direction toward the interior.

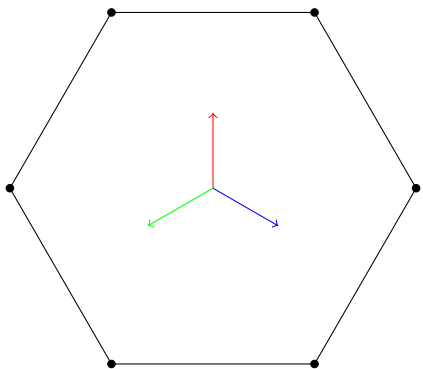

FIGURE 7. A hexagon as the space of game states, with arrows representing the directions available for the random walk.

This hexagon model is created to address a potential issue of the triangle model: as the random walk approaches a vertex, the probability of being ejected (hitting an edge) swings wildly between the two losing players, which may be undesirable. The hexagon model, if ejection is implemented, allows for the ejection zones to be separated by a nonzero distance, avoiding this issue.

Another way to model the three-player game can be represented as a complete ternary tree with some depth $d$, corresponding to the fineness of the time steps. Each leaf vertex is assigned Red, Green, or Blue (corresponding to the players), in such a way that nodes sharing the same parent must have different colors, see Figure 8.

The game proceeds by tree traversal: starting at the root, a particle moves a level down the tree during each time step, reaching each child with probability $1/3$. During each step, a lucky player is chosen uniformly at random, and some subset of the descendant leaves are converted to the color corresponding to that player. The game ends when a leaf is reached, with the color of the leaf corresponding to the winner.

5. CONCLUDING REMARKS; LIMITATIONS AND FUTURE WORK

Furthering the work of [KS] and [BM] on a stochastic game inspired by backgammon, we determined decision thresholds for more general variants of game where the number of doubles and stake multipliers may vary. We also introduced a three-player generalization of the game and proved basic results about game behavior. We saw that in general a double (or its generalization) should be proposed when the opponent(s) is indifferent between accepting and rejecting the proposal; this principle of indifference allowed the explicit computation of decision thresholds.

Our results are limited in their applicability to an actual game of backgammon.

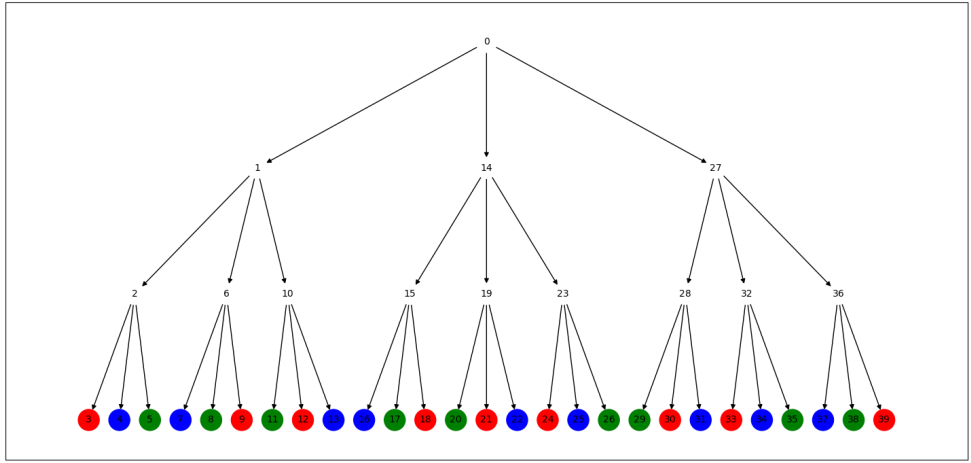

FIGURE 8. Initial state of the ternary tree model.

- Our random walk models assume that the walk is continuous or discretized in a way such that only a small distance may be traversed in the span of a few turns. The game evaluation of an actual game of backgammon may change much more erratically: one move, especially in the endgame, could result in a wild swing in the evaluation.
- We assume that players (i) have the sole objective of maximizing expected payoff and (ii) play perfectly toward that objective, which might not be the case in an actual game. A player close to winning a tournament may choose to play conservatively to maintain their ranking (an objective not considered by our models), while someone who fared more poorly may be more open to risks for the possibility of placing in the tournament. Furthermore, against a less skilled opponent, it may be sensible to continue the game instead of forfeiting when the threshold is reached, in hopes of a future blunder by the opponent.
- The expected payoff, in practice, is relevant only over the span of many games, enough so that the law of large numbers is effective. If only a small number of games is played, a player may fare better with a lower-risk, lower-reward strategy.
- As there are no pieces to move in our model, there is no notion of "going first"; the initial game state of $0.5$ reflects the assumption that the game starts off balanced, and neither player has an advantage over the other. In an actual game of backgammon, however, the first mover has a slight advantage. In comparison, having the first move in Chess, although unproven, is hypothesized to be advantageous due to the initiative, and having the first move in Giveaway Chess has been shown to be winning.
- Our thresholds are computed using the game evaluation, which in our model is universal and visible to every player. In an actual game of backgammon, player do not have access to this information and must make an educated guess based on the game state; depending on personal factors such as playing style and skill level, different players may arrive at different evaluations of the same position.

Our analysis points naturally to some directions of further work.

- Our current random walk models do not account for the wild swings in probability that may occur in an actual backgammon game. One could introduce a more sophisticated model where, at each time step, both the direction and magnitude of the step taken are randomized.

- In the three-player variant based on random walk in a triangle, we were only able to determine the decision regions for one double only. The difficulty, as pointed out in Remark 4.5, was determining the probability of a 2-dimensional random walk hitting one region before another. If further progress can be made toward this sub-problem, it would be feasible to determine the decision regions for more complex doubling rules. It would be interesting to see if the decision regions have a limiting shape in the case of infinitely many doubles.

    Furthermore, it would be interesting to replace the 6-directional discrete random walk with a continuous random walk, such as Brownian motion.

- It would be interesting to consider variants with more than three players. The 3-player doubling rules in Section 4 generalize immediately to more players, but it is unclear what polytope to choose for modeling the game state. A higher-order two-dimensional polygon runs into a problem that the evaluations of players in close proximity are correlated, so perhaps a higher-dimensional simplex would be the most natural choice.

- In the discrete random walk in $\triangle ABC$ where the equilateral triangle tiling has edges subdivided by $n$, each vertex of the tiling (and thus each game state and evaluation) can be written as $(a/n)A + (b/n)B + (c/n)C$, where $a+b+c = n$ are integers. The triple $(a/n, b/n, c/n)$ can also be thought of as the type (empirical distribution) of a string with length $n$ over the alphabet $\{A, B, C\}$. This connection may be worth further investigation, to give a combinatorial interpretation to the random walk, in addition to a geometric one. Furthermore, such a perspective may make information-theoretic tools applicable.

## Appendix A. Proofs of auxiliary results

*Proof of Theorem 1.1.* Since Player B plays perfectly, Player A's expected payoff is never greater than the current stake of the game. Thus, whenever Player B is willing to fold, Player A should double, since doing so results in receiving the current stake of the game, which is the greatest amount of payoff possible. Thus the optimal point to double should be at least as early as any point when Player B is willing to fold, so $\alpha \leq \beta$.

Suppose that Player A doubles at some point $p$ strictly less than $\beta$, in which case Player B would accept. Let $q$ be such that $p < q < \beta$; after the doubling and acceptance, Player A's state may reach $q$, or it may not. If Player A's state does reach $q$, Player A could have waited till $q$ to propose the double, since the result would have been the same. If Player A's state does not reach $q$, then Player A loses, and loses double the stake. Thus, doubling at $p$ is inferior to doubling at $q$, which implies the optimal time to double cannot be less than $\beta$. Thus $\alpha \geq \beta$, which together with $\alpha \leq \beta$ implies $\alpha = \beta$. $\square$

*Proof of Lemma 3.3.* Suppose Player B raises. Player A, at evaluation $a_k$ ($k \geq 2$), is indifferent between accepting and folding. If Player A folds, then their expected payoff is $-x$. If Player A accepts, either they lose and receive $-y$, or their evaluation rises to $d_{k-1}$ and they offer the $(k-1)^{\text{st}}$ raise, at which point the expected payoff is $yx$. The probability of reaching 0 before $d_{k-1}$ is $(d_{k-1} - a_k)/d_{k-1}$, while the probability of reaching $d_{k-1}$ before 0 is $a_k/d_{k-1}$. Therefore, the expected payoff for Player A, upon acceptance, is

$$-y \left( \frac{d_{k-1} - a_k}{d_{k-1}} \right) + yx \left( \frac{a_k}{d_{k-1}} \right). \tag{A.1}$$

By Theorem 1.3, the threshold occurs when Player A's expected payoff is the same whether the "double" is accepted or not. Setting (A.1) equal to $-x$, and noting that $a_k = 1 - d_k$, yields

$$d_k = 1 - d_{k-1} \left( \frac{y-x}{y(x+1)} \right). \tag{A.2}$$

To find the limiting threshold as the number of raises becomes large, it suffices to solve for the value where the recurrence stabilizes:

$$p = 1 - p\left(\frac{y-x}{y+yx}\right)$$

$$\Longrightarrow p = \frac{1}{1 + \frac{y-x}{y+yx}} = \frac{y+yx}{2y+yx-x}. \tag{A.3}$$

To find a closed-form for the threshold given arbitrarily many offers remaining, one can try to track the error term $\epsilon_k := d_k - p$, and obtain

$$\begin{aligned}
\epsilon_k &= 1 - d_{k-1}\left(\frac{y-x}{y+yx}\right) - \frac{y+yx}{2y+yx-x} \\
&= \frac{y-x}{2y+yx-x} - \frac{y-x}{y+yx}\left(d_{k-1} - \frac{y+yx}{2y+yx-x} + \frac{y+yx}{2y+yx-x}\right) \\
&= \frac{y-x}{2y+yx-x} - \frac{y-x}{2y+yx-x} - \left(\frac{y-x}{y+yx}\right)\epsilon_{k-1} \\
&= \left(-\frac{y-x}{y+yx}\right)\epsilon_{k-1}. \tag{A.4}
\end{aligned}$$

Therefore the error terms satisfies $\epsilon_k = c\left(-\frac{y-x}{y+yx}\right)^k$ for some constant $c$. From the initial condition $\epsilon_1 + p = d_1 = 1/2 + x/(2y)$, given in Lemma 3.2, we can determine $c$ and obtain that

$$d_k = \frac{y(x+1)}{2y+x(y-1)} + \frac{x(x+1)(y-1)}{2(2y+xy-x)}\left(-\frac{y-x}{y(x+1)}\right)^k. \tag{A.5}$$

$\square$

*Proof of Lemma 3.4.* Suppose Player B reduces. Player A, at evaluation $a_k$ ($k \geq 2$), is indifferent between accepting and folding. If Player A folds, then their expected payoff is $-x$. If Player A accepts, either they win and receive $y$, or their evaluation drops to $d_{k-1}$ and they offer the $(k-1)^{st}$ reduction, at which point the expected payoff is $yx$. The probability of reaching 1 before $d_{k-1}$ is $(a_k - d_{k-1})/(1 - d_{k-1})$, while the probability of reaching $d_{k-1}$ before 1 is $(1 - a_k)/(1 - d_{k-1})$. Therefore, Player A's expected payoff upon acceptance is

$$y\left(\frac{a_k - d_{k-1}}{1 - d_{k-1}}\right) + yx\left(\frac{1 - a_k}{1 - d_{k-1}}\right). \tag{A.6}$$

By Theorem 1.3, the threshold occurs when Player A's expected payoff is the same whether the reduction is accepted or not. Setting (A.6) equal to $-x$, and noting that $a_k = 1 - d_k$, yields

$$d_k = \frac{y+x}{y(1-x)}(1 - d_{k-1}). \tag{A.7}$$

To find the limiting threshold as the number of reductions becomes large, it suffices to solve for the value where the recurrence stabilizes:

$$p = (1-p)\left(\frac{y+x}{y(1-x)}\right)$$

$$\Longrightarrow p = \frac{y+x}{2y-yx+x}. \tag{A.8}$$

To find a closed-form for the threshold given arbitrarily many offers remaining, one can perform a error analysis on $\epsilon_k := d_k - p$, similarly to the proof of Lemma 3.3, and obtain

$$\epsilon_k = \left(-\frac{y+x}{y(1-x)}\right)\epsilon_{k-1}. \tag{A.9}$$

The constant coefficient can be determined from the $k = 1$ case given by Lemma 3.2, and from there it follows that

$$d_k = \frac{y+x}{2y-yx+x} + \frac{(-x)(1-x)(1-y)}{2(2y-yx+x)}\left(-\frac{y+x}{y(1-x)}\right)^k. \tag{A.10}$$

$\square$

APPENDIX B. MODELLING 3-PLAYER GAME WITH A 2-DIMENSIONAL RANDOM WALK

Here we list the code used to simulate the 3-player game with a 2D random walk.

The first function simulates one round of the game using a 2D random walk. We input a desired radius and starting coordinate, normally the origin, and it will give us a winner who will either be Red, Green, or Blue, by performing this walk. The base for this walk comes from [Han] and we adapt their walk from movement in four directions to moving in six directions according to a triangular grid.

```python
import matplotlib.pyplot as plt
import random
import numpy as np

# Define our walk to simulate one round of Backgammon
# Parameters are radius and initial x/y coords
# We define the radius as the maximum distance you could travel in a straight line from the origin
# First we create an array to store all the directions our walk can travel in
# We allow our walk to move in 6 directions

directions = [(0,1),(3**(1/2)/2,-1/2),(-3**(1/2)/2,-1/2),(0,-1),(-3**(1/2)/2,1/2),(3**(1/2)/2,1/2)]

# We can now define our first function

def Round(r,x,y):

    # Create arrays to store the walk

    directions_x = []
    directions_y = []

    # Define our initial coordinates

    current_step_x = x
    current_step_y = y

    # Start a while loop, looping over the max number of steps n

    while True:

        # Insert our win condition for if it escapes the circle
        # First check if the walk has escaped the circle
```

```python
33          # We create a boundary just before the edge to account for floating numbers
34          # First deal with the three points where each player wins
35  
36          if current_step_y >= current_step_x/(3**(1/2)) + r and current_step_y >= r -
     current_step_x/(3**(1/2)):
37              return 'B'
38  
39          elif current_step_y <= -current_step_x/(3**(1/2)) - r and current_step_x <= -r
     *3**(1/2)/2:
40              return 'G'
41  
42          elif current_step_y <= current_step_x/(3**(1/2)) - r and current_step_x >= r
     *3**(1/2)/2:
43              return 'R'
44  
45          # We now consider the three edges/boundaries
46          # Each one of these is considered a two player game between the two winners of
     the neighbouring edges
47          # We take the x value of our final coord and calculate how far along the edge
     it is
48          # We use this to find the probability of each player winning
49  
50          elif current_step_y <= -r/2 + 1/4 and -r*3**(1/2)/2 <= current_step_x <= r
     *3**(1/2)/2:
51  
52              p = (current_step_x + r*3**(1/2)/2)/(r*3**(1/2))
53              result = random.choices(['R','G'], weights = (p,1-p))
54              return result[0]
55  
56          elif current_step_y >= r + current_step_x*3**(1/2) -1/4 and current_step_y >=
     -current_step_x/(3**(1/2)) - r and current_step_y <= r - current_step_x/(3**(1/2)):
57  
58              p = (current_step_x + r*3**(1/2)/2)/(r*3**(1/2)/2)
59              result = random.choices(['B','G'], weights = (p,1-p))
60              return result[0]
61  
62          elif current_step_y >= r - 3**(1/2)*current_step_x - 1/4 and current_step_y >=
      current_step_x/(3**(1/2)) - r and current_step_y <= current_step_x/(3**(1/2)) + r:
63  
64              p = (current_step_x)/(r*3**(1/2)/2)
65              result =  random.choices(['R','B'], weights = (p,1-p))
66              return result[0]
67  
68          # If the code reaches this point, it has not yet escaped the triangle
69          # Define the walk
70  
71          else:
72  
73              # Randomly pick a direction from our 6 choices
74  
75              chosen_step = random.choice(directions)
76  
77              # Add the vector direction to our current coordinate
```

```
78
79              current_step_x += chosen_step[0]
80              current_step_y += chosen_step[1]
81
82              # Add the new coordinate to our random walk
83
84              directions_x.append(current_step_x)
85              directions_y.append(current_step_y)
86
87      # If the code reaches this point, it never leaves the circle after all steps have moved
88      # Return N for no winner
89
90      return 'N'
```

Figure 9 shows an example output from this round function played ten times. We can see Blue wins five times, Green wins four times, and Red wins once.

```
B
B
G
G
R
B
G
B
G
B
```

FIGURE 9. Output from Round function played ten times, with initial conditions $(10, 0, 0)$

Now we can simulate a round of backgammon, we can use this to calculate the probability of a player winning from any point on the lattice. We can do this by simulating many games and calculating what proportion of times each player wins from that position, varying the initial conditions as allowed by our round function above.

```
1 # This simulates one round
2 # Define a new function which simulates several rounds
3 # Inputs are radius, no of games, and initial x/y coords
4
5 def Backgammon_Triangle(r,g,x,y):
6
```

```
7       # Initialise arrays to store the number of wins for each player and none
8
9       green = []
10      blue = []
11      red = []
12      none = []
13
14      # Initialise our loop variable
15
16      i = 0
17
18      # Loop over the number of games
19
20      while i < g:
21
22          # For each i, we simulate a game and find the winner
23
24          result = Round(r,x,y)
25
26          # Check the winner of each round
27          # Whoever wins, we add one (element) to their array
28          # Then add one to our loop variable to continue the loop
29
30          if result == 'G':
31              green.append(1)
32              i += 1
33          elif result == 'B':
34              blue.append(1)
35              i += 1
36          elif result == 'R':
37              red.append(1)
38              i += 1
39          elif result == 'N':
40              none.append(1)
41              i += 1
42
43      # Create an array which stores the experimental probability of each player of
        winning
44      # We do this by taking the sum of each array and dividing it by the number of
        games
45
46      final_result = [sum(green)/g,sum(blue)/g,sum(red)/g]
47
48      # Return this array
49
50      return final_result
```

Figure 10 an example output from the Backgammon_Triangle function played three times with various numbers of game played. We started these functions at the origin, and so would expect the result to be approximately $1/3$ for each. We can see that the more games we play, the more accurate our function is, as we would expect.

We now have the probability of each player winning from any location on the lattice. We can now introduce the notion of doubling as discussed in the preliminary paper. To do this, we have to make some

```
[0.45, 0.32, 0.23]
[0.338, 0.337, 0.325]
[0.3378, 0.3267, 0.3355]
```

FIGURE 10. Output from Backgammon_Triangle function called three times, with initial conditions $(10, g, 0, 0)$, with $g = 10, 100, 1000$ in that order.

decisions on how we introduce doubling. In this version, we decided that both players needed to accept the double for the game to continue, and if even one rejects the game ends and the offer-er wins. We also decided that anyone could double, as long as they weren't the last player to double. In this final part, we create a function which considers two probabilities; a minimum probability for a player to offer a double, and a minimum probability for a player to accept a double. We can see an example output in Figure 11

```python
# We now want to introduce a stake and doubling system
# We define a function whose input is radius, the min prob to offer a double, and the
    min prob to accept a double

def Doubling(r,q_offer,q_reject):

    # Define our initial stake as 1 and introduce a variable to keep track of the most
     recent doubler

    stake = 1
    doubler = None

    # Create arrays to store the walk

    directions_xx = [0]
    directions_yy = [0]

    # Define our initial coordinates

    current_step_xx = 0 #current x-coordinate
    current_step_yy = 0 #current y-coordinate

    # Set up plots

    figure, axes = plt.subplots()

    hexagon_x1 = [0,r*3**(1/2)/2,r*-3**(1/2)/2,0]
    hexagon_y1 = [r,r*-1/2,r*-1/2,r]
    axes.plot(hexagon_x1,hexagon_y1,color='k')

    hexagon_x2 = [0,r*3**(1/2)/2 -1/2,r*-3**(1/2)/2 +1/2,0]
    hexagon_y2 = [r -1/2,r*-1/2 +1/4,r*-1/2 +1/4,r -1/2]
    axes.plot(hexagon_x2,hexagon_y2,color='r')

```

```
33        for i in range(r):
34            hexagon_x_r = [0,i*3**(1/2)/2,i*-3**(1/2)/2,0]
35            hexagon_y_r = [i,i*-1/2,i*-1/2,i]
36            axes.plot(hexagon_x_r,hexagon_y_r,color='0.8')
37
38        # Start a while loop, looping over the number of steps
39
40        while True:
41
42            # We reinnsert our win condition for if it escapes the circle
43            # This is the same as in our Round function
44            # We format the return to include the current stake
45            # We also print the probability if it lands on an inbetween boundary
46
47            if current_step_yy >= current_step_xx/(3**(1/2)) + r and current_step_yy >= r
        - current_step_xx/(3**(1/2)):
48                return print('B wins {}'.format(stake))
49
50            elif current_step_yy <= -current_step_xx/(3**(1/2)) - r and current_step_xx <=
         -r*3**(1/2)/2:
51                return print('G wins {}'.format(stake))
52
53            elif current_step_yy <= current_step_xx/(3**(1/2)) - r and current_step_xx >=
        r*3**(1/2)/2:
54                return print('R wins {}'.format(stake))
55
56            elif current_step_yy <= -r/2 + 1/4 and -r*3**(1/2)/2 <= current_step_xx <= r
        *3**(1/2)/2:
57
58                p = (current_step_xx + r*3**(1/2)/2)/(r*3**(1/2))
59                print('R has a ' + str(p*100) + '% chance of winning')
60                result = random.choices(['R','G'], cum_weights = (p,1), k = 1)
61                return print('{} wins {}'.format(result[0],stake))
62
63            elif current_step_yy >= r + current_step_xx*3**(1/2) - 1/4 and current_step_yy
         >= -current_step_xx/(3**(1/2)) - r and current_step_yy <= r - current_step_xx
        /(3**(1/2)):
64
65                p = (current_step_xx + r*3**(1/2)/2)/(r*3**(1/2)/2)
66                print('B has a ' + str(p*100) + '% chance of winning')
67                result = random.choices(['B','G'], cum_weights = (p,1), k = 1)
68                return print('{} wins {}'.format(result[0],stake))
69
70            elif current_step_yy >= r - 3**(1/2)*current_step_xx - 1/4 and current_step_yy
         >= current_step_xx/(3**(1/2)) - r and current_step_yy <= current_step_xx/(3**(1/2)
        ) + r:
71
72                p = (current_step_xx)/(r*3**(1/2)/2)
73                print('R has a ' + str(p*100) + '% chance of winning')
74                result =  random.choices(['R','B'], cum_weights = (p,1) , k = 1)
75                return print('{} wins {}'.format(result[0],stake))
76
77        # If it reaches this point, the walk has not yet escaped the circle
```

```python
        else:
            
            # We calculate the probability of each player winning from the current point
            # We do this by simulating 10000 games from that position and count the winners
            # We use our Backgammon_Hexagon function to do this
            # We then print the current probability, stake, and last doubler
            
            prob = Backgammon_Triangle(r,10000,current_step_xx,current_step_yy)
            print(prob,stake,doubler)
            
            # For player G, we check if they are above the min prob to offer a double and whether they are not the doubler
            
            if prob[0]>q_offer and doubler != 'G':
                
                # If both checks are passed, we check if either of the other plays are below the min prob to accept the double
                
                if prob[1]<q_reject or prob[2]<q_reject:
                    
                    # If either reject, the game ends and the doubler wins
                    # We plot the walk and return the winner and the stake they won
                    
                    axes.plot(directions_xx,directions_yy,color="k")
                    plt.show()
                    return print('G wins {}'.format(stake))
                
                else:
                    
                    # If both accept, the stake doubles and we update who the current doubler is
                    
                    stake = 2*stake
                    doubler = 'G'
            
            # We do the same for player B
            
            
            elif prob[1]>q_offer and doubler != 'B':
                
                if prob[0]<q_reject or prob[2]<q_reject:
                    
                    axes.plot(directions_xx,directions_yy,color="k") #make plot using lists of x-, y-coordinates
                    plt.show()
                    return print('B wins {}'.format(stake))
                
                else:
                    
                    stake = 2*stake
```

```
125                    doubler = 'B'
126
127            # We do the same for player R
128
129            elif prob[2]>q_offer and doubler != 'R':
130
131                if prob[0]<q_reject or prob[1]<q_reject:
132
133                    axes.plot(directions_xx,directions_yy,color="k") #make plot using
     lists of x-, y-coordinates
134                    plt.show()
135                    return print('R wins {}'.format(stake))
136
137                else:
138
139                    stake = 2*stake
140                    doubler = 'R'
141
142        # If it reaches this point, the walk continues
143        # Randomly pick a direction from our 6 choices
144
145        chosen_step = random.choice(directions)
146
147        # Add the vector direction to our current coordinate
148
149        current_step_xx += chosen_step[0]
150        current_step_yy += chosen_step[1]
151
152        # Add the new coordinate to our random walk
153
154        directions_xx.append(current_step_xx)
155        directions_yy.append(current_step_yy)
156
157        # Plot our walk so far
158
159        axes.plot(directions_xx,directions_yy,color="k")
```

In our final function, we explore how the probability of a player winning compares to the area of their triangle. We first introduce a function which calculates the area of a triangle given three coordinates. The final function then uses this and an adapted version of the doubling function to compare these results. See Figure 12 for an example output.

```
1  def Triangle_Area(x1,y1,x2,y2,x3,y3):
2      area = 1/2*abs(x1*(y2-y3)+x2*(y3-y1)+x3*(y1-y2))
3      return area
4
5  def Prob_Area(r):
6
7      # Create arrays to store the walk
8
9      directions_xx = [0]
10     directions_yy = [0]
11
12     # Define our initial coordinates
```

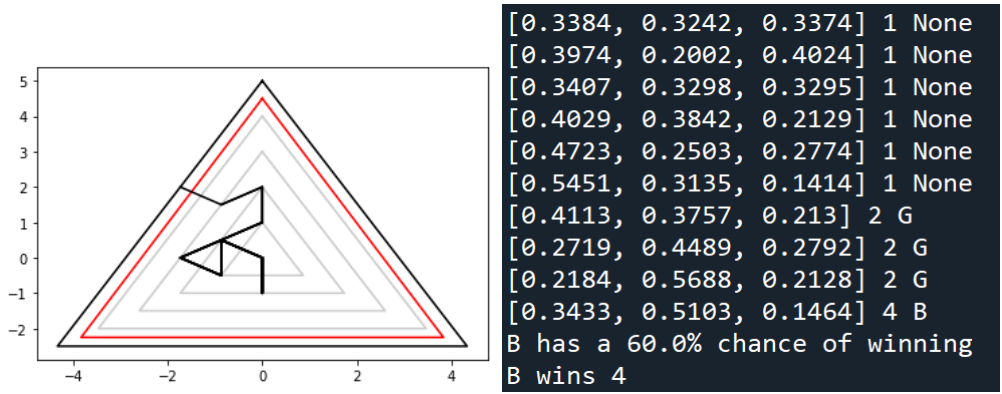

FIGURE 11. Output from Doubling function called three times, with initial conditions $(5, 0.5, 0.1)$, showing the random walk and the probability of each player winning at each point visited alongside the current stake and the current doubler.

```
13
14      current_step_xx = 0 #current x-coordinate
15      current_step_yy = 0 #current y-coordinate
16
17      # Set up plots
18
19      figure, axes = plt.subplots()
20
21      hexagon_x1 = [0,r*3**(1/2)/2,r*-3**(1/2)/2,0]
22      hexagon_y1 = [r,r*-1/2,r*-1/2,r]
23      axes.plot(hexagon_x1,hexagon_y1,color='k')
24
25      hexagon_x2 = [0,r*3**(1/2)/2 -1/2,r*-3**(1/2)/2 +1/2,0]
26      hexagon_y2 = [r -1/2,r*-1/2 +1/4,r*-1/2 +1/4,r -1/2]
27      axes.plot(hexagon_x2,hexagon_y2,color='r')
28
29      for i in range(r):
30          hexagon_x_r = [0,i*3**(1/2)/2,i*-3**(1/2)/2,0]
31          hexagon_y_r = [i,i*-1/2,i*-1/2,i]
32          axes.plot(hexagon_x_r,hexagon_y_r,color='0.8')
33
34      # Start a while loop, looping over the number of steps
35
36      while True:
37
38          # We reinnsert our win condition for if it escapes the circle
39          # This is the same as in our Round function
40          # We format the return to include the current stake
41          # We also print the probability if it lands on an inbetween boundary
42
43          if current_step_yy >= current_step_xx/(3**(1/2)) + r and current_step_yy >= r - current_step_xx/(3**(1/2)):
44              return print('B wins')
45
```

```python
        elif current_step_yy <= -current_step_xx/(3**(1/2)) - r and current_step_xx <= -r*3**(1/2)/2:
            return print('G wins')

        elif current_step_yy <= current_step_xx/(3**(1/2)) - r and current_step_xx >= r*3**(1/2)/2:
            return print('R wins')

        elif current_step_yy <= -r/2 + 1/4 and -r*3**(1/2)/2 <= current_step_xx <= r*3**(1/2)/2:
            p = (current_step_xx + r*3**(1/2)/2)/(r*3**(1/2))
            print('R has a ' + str(p*100) + '% chance of winning')
            result = random.choices(['R','G'], cum_weights = (p,1), k = 1)
            return print('{} wins'.format(result[0]))

        elif current_step_yy >= r + current_step_xx*3**(1/2) - 1/4 and current_step_yy >= -current_step_xx/(3**(1/2)) - r and current_step_yy <= r - current_step_xx/(3**(1/2)):
            p = (current_step_xx + r*3**(1/2)/2)/(r*3**(1/2)/2)
            print('B has a ' + str(p*100) + '% chance of winning')
            result = random.choices(['B','G'], cum_weights = (p,1), k = 1)
            return print('{} wins'.format(result[0]))

        elif current_step_yy >= r - 3**(1/2)*current_step_xx - 1/4 and current_step_yy >= current_step_xx/(3**(1/2)) - r and current_step_yy <= current_step_xx/(3**(1/2)) + r:
            p = (current_step_xx)/(r*3**(1/2)/2)
            print('R has a ' + str(p*100) + '% chance of winning')
            result =  random.choices(['R','B'], cum_weights = (p,1) , k = 1)
            return print('{} wins'.format(result[0]))

        # If it reaches this point, the walk has not yet escaped the circle

        else:

            # We calculate the probability of each player winning from the current point
            # We do this by simulating 10000 games from that position and count the winners
            # We use our Backgammon_Triangle function to do this
            # We then calculate the area of the triangle of each player
            # By dividing by the total area of the triangle, we can compare this to the probability of each player winning (assuming they are equivalent)

            prob = Backgammon_Triangle(r,10000,current_step_xx,current_step_yy)

            area_R = Triangle_Area(current_step_xx,current_step_yy,0,r,-r*3**(1/2)/2,-r*1/2)
            area_G = Triangle_Area(current_step_xx,current_step_yy,0,r,r*3**(1/2)/2,-r*1/2)
```

```
87              area_B = Triangle_Area(current_step_xx,current_step_yy,r*3**(1/2)/2,-r
        *1/2,-r*3**(1/2)/2,-r*1/2)
88              area = Triangle_Area(0,r,r*3**(1/2)/2,-r*1/2,-r*3**(1/2)/2,-r*1/2)
89
90              relative_area = [area_G/area,area_B/area,area_R/area]
91              print(np.subtract(prob,relative_area))
92
93          # If it reaches this point, the walk continues
94          # Randomly pick a direction from our 6 choices
95
96          chosen_step = random.choice(directions)
97
98          # Add the vector direction to our current coordinate
99
100         current_step_xx += chosen_step[0]
101         current_step_yy += chosen_step[1]
102
103         # Add the new coordinate to our random walk
104
105         directions_xx.append(current_step_xx)
106         directions_yy.append(current_step_yy)
107
108         # Plot our walk so far
109
110         axes.plot(directions_xx,directions_yy,color="k")
```

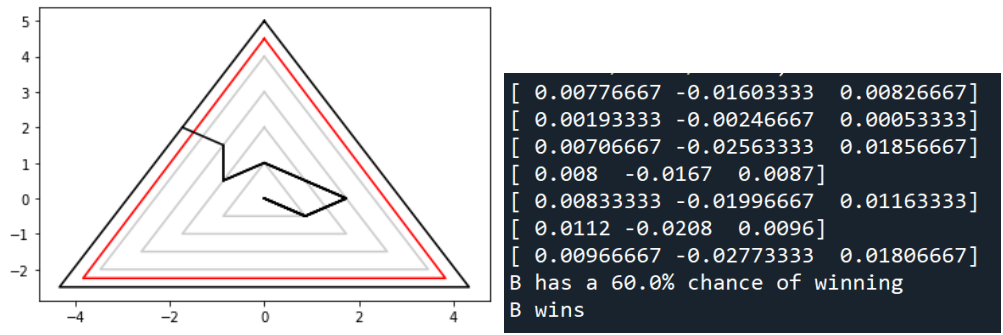

FIGURE 12. Output from Prob_Area function called three times, with initial conditions (5), showing the random walk and the difference between the probability of each player winning at that point and the area of their triangle.

## REFERENCES

[BM]    J. Blath and P. Mörters, *Backgammon, doubling the stakes, and Brownian motion*, Plus Magazine (University of Cambridge), 01 May 2001. https://plus.maths.org/content/backgammon-doubling-stakes-and-brownian-motion.

[Don]    M. D. Donsker, *An invariance principle for certain probability limit theorems*, Memoirs of the American Mathematical Society, 6 (1951), pp. 1–10.

[Han]    S. Hansen, *Constructing random walks with Python*, From Code to Clinic, Retrieved from https://fromcodetoclinic.wordpress.com/2017/01/25/random-walks-w-python/

[KS]    E. B. Keeler and J. Spencer, *Optimal Doubling in Backgammon*, Operations Research, Vol. 23, No. 6 (1975), pp. 1063-1071.

*Email address*: jeffreyju1899@gmail.com

DEPARTMENT OF MATHEMATICS, UNIVERSITY OF SOUTHERN CALIFORNIA, LOS ANGELES, CA 90007

*Email address*: dleifer@umich.edu

DEPARTMENT OF MATHEMATICS, UNIVERSITY OF MICHIGAN, ANN ARBOR, MI 48109

*Email address*: sjm1@williams.edu

DEPARTMENT OF MATHEMATICS AND STATISTICS, WILLIAMS COLLEGE, WILLIAMSTOWN, MA 01267

*Email address*: soorajap@iitb.ac.in

DEPARTMENT OF MATHEMATICS, INDIAN INSTITUTE OF TECHNOLOGY, BOMBAY

*Email address*: csun1248@gmail.com

DEPARTMENT OF MATHEMATICS, COLUMBIA UNIVERSITY, NEW YORK, NY 10027

*Email address*: lt4@williams.edu

DEPARTMENT OF MATHEMATICS AND STATISTICS, WILLIAMS COLLEGE, WILLIAMSTOWN, MA 01267

*Email address*: bentocher42@gmail.com

MATHEMATICS INSTITUTE, UNIVERSITY OF WARWICK, COVENTRY, CV4 7AL

*Email address*: kpw1@williams.edu

DEPARTMENT OF MATHEMATICS AND STATISTICS, WILLIAMS COLLEGE, WILLIAMSTOWN, MA 01267